\documentclass[12pt]{amsart}
\usepackage[all]{xy}
\usepackage{amssymb}
\calclayout
\newtheorem{dummy}{anything}[section] 
\newtheorem{theorem}[dummy]{Theorem}
\newtheorem*{thma}{Theorem A}
\newtheorem*{thmb}{Theorem B}
\newtheorem{lemma}[dummy]{Lemma} 
\newtheorem{proposition}[dummy]{Proposition} 
\newtheorem{corollary}[dummy]{Corollary}
 
\theoremstyle{definition}%%Change Theoremstyle
\newtheorem{definition}[dummy]{Definition}

 \newtheorem{remark}[dummy]{Remark}

\newcommand
{\eqncount}{\setcounter{equation}{\value{dummy}}%
\addtocounter{dummy}{1}}
\newcommand{\cH}{\mathcal H}

\newcommand{\cS}{\mathcal S}

\newcommand{\cE}{\mathcal E}
\newcommand{\cT}{\mathcal T}

\newcommand{\bZ}{\mathbf Z}

\newcommand{\bC}{\mathbb C}

\newcommand{\cy}[1]{\bZ/{#1}}

\newcommand{\bd}{\partial}
\newcommand{\vv}{\, | \,}
\newcommand{\trf}{tr{\hskip -1.8truept}f}

\newcommand{\mmatrix}[4]{\left (\vcenter
{\xymatrix@C-2pc@R-2pc{#1&#2\\#3&#4} }
\right )}

\newcommand{\disjointunion}{\hbox{$\perp\hskip -4pt\perp$}}
\DeclareMathOperator{\Hom}{Hom}

\DeclareMathOperator{\rank}{rank}
\DeclareMathOperator{\Sharp}{\sharp}
\DeclareMathOperator{\Image}{Im}

\DeclareMathOperator{\hepta}{Aut}
\newcommand{\he}[1]{\hepta(#1)}
\newcommand{\hept}[1]{\hepta_{\bullet}(#1)}
\newcommand{\heq}[1]{\cE (#1)}
\newcommand{\heqpt}[1]{\cE_{\bullet} (#1)}
\newcommand{\htildeM}{\widetilde \cH (M)}
\newcommand{\htildeB}{\widetilde \cH (B)}

\newcommand{\hM}{\cH (M)}
\newcommand{\hB}{\cH (B)}
\DeclareMathOperator{\Isom}{Isom}

\newcommand{\quadtypeM}{[\pi_1, \pi_2, k_M, s_M]}
\newcommand{\Ospin}{\Omega^{Spin}}
\newcommand{\rOspin}{\widehat\Omega^{Spin}}

\newcommand{\Kw}{KH_2(M;\cy 2)}
\newcommand{\Bw}{B\langle w_2\rangle}
\newcommand{\Mw}{M\negthinspace\langle w_2\rangle}
\newcommand{\whept}[1]{\hepta_{\bullet}(#1,w_2)}
\newcommand{\whtildeM}{\widetilde \cH (M,w_2)}
\newcommand{\whtildeB}{\widetilde \cH (B,w_2)}

\begin{document}

\title[Homotopy Equivalences]
{Homotopy Self-Equivalences of $4$-manifolds}
\author{Ian Hambleton and Matthias Kreck} 
\thanks{\hskip -11pt Research partially supported by 
NSERC Research Grant A4000. The authors thank the Max Planck Institut f\"ur Mathematik in Bonn, Germany, for hospitality and
support.
}
\address{Department of Mathematics \& Statistics
 \newline\indent
McMaster University
 \newline\indent
Hamilton, ON  L8S 4K1, Canada}
\email{ian{@}math.mcmaster.ca}
\address{Mathematisches Institut
 \newline\indent
Universit\"at Heidelberg
 \newline\indent
D-69120 Heidelberg, Germany}
\email{kreck{@}mathi.uni-heidelberg.de}

\date{Nov. 18, 2003}
\begin{abstract}\noindent
We establish a braid of interlocking exact sequences containing the group of homotopy self-equivalences
of a smooth or topological $4$-manifold. The braid is computed 
for manifolds whose fundamental group  is finite of odd order.
\end{abstract}
\maketitle
\section{Introduction}
Let $M^4$ be a closed, oriented, smooth or topological $4$-manifold.
We wish to study the group $\he{M}$ of homotopy
classes of 
 homotopy self-equivalences
$f\colon M \to M$,
using techniques from surgery and bordism
theory. We will always assume that $M$ is connected.
Here is an overview of our results, starting with
 an informal description of some related objects. 

The group $\hM$ consists of 
oriented $h$-cobordisms $W^5$ from $M$ to $M$, under the
equivalence relation induced by $h$-cobordism relative to the
boundary. The orientation of $W$ induces  opposite
orientations on the two boundary components $M$.
An $h$-cobordism gives a homotopy self-equivalence of $M$,
and we get a homomorphism $\hM\to \he{M}$. 

Let $B = B(M)$ denote the $2$-type of $M$. It
is a fibration over $K(\pi_1(M),1)$, with fibre
$\pi_2(M)$ determined by a $k$-invariant
$k_M\in H^3(\pi_1;\pi_2)$, obtained from $M$ by
attaching cells of dimension $\geq 4$ to kill the homotopy groups in
dimensions $\geq 3$. The natural map
$c\colon M \to B$ is $3$-connected, and we refer
to this as the classifying map of $M$. There is an
induced homomorphism $\he M \to \he B$,
the group of homotopy
classes of homotopy self-equivalences of $B$, by 
obstruction theory and the
naturality of the construction.
If $M$ is a spin manifold, 
we will also be using the  smooth (or topological) 
bordism groups $\Ospin_n(B)$. 
By imposing the requirement that the reference maps
to $M$ must have degree zero, we obtain modified bordism
groups $\rOspin_4(M)$ and $\rOspin_5(B,M)$.
When
$w_2(M)
\neq 0$, we will use the appropriate bordism groups
of the \emph{normal $2$-type}.
See \cite{stong1}, \cite{kreck3}, \cite{kreck2} for this theory.

A variation of $\hM$, denoted $\htildeM$, will
also be useful. This is the group of oriented
bordisms $(W,\bd_{-}W, \bd_{+} W)$ with
$\bd_{\pm}W = M$, equipped with a map
$F\colon W \to M$.
We require  the restrictions $F|_{\bd_{\pm}W}$
 to the boundary components
to be   homotopy equivalences (and the identity on
 the component $\bd_{-} W$). The equivalence
relation on these objects is induced by
bordism (extending the map to $M$) relative to the boundary
(see Section 2.2  for the details).

Our strategy is to compare $\he M$ to these other groups by means of
various interlocking exact sequences. For technical reasons, we will
restrict ourselves to homotopy self-equivalences preserving both
the given orientation on $M$ and a 
fixed base-point $x_0 \in M$. 
Let $\hept M$ denote the group of
homotopy classes of such homotopy self-equivalences.
We will also define
 ``pointed" versions of the other objects, including
the space $\heqpt{B}$ of base-point preserving homotopy
equivalences of $B$, and 
 the group $\hept{B}=\pi_0(\heqpt{B})$.
Our main qualitative result in the spin case, Theorem \ref{mainspin}, is expressed in a commutative braid
\vskip .3cm
$$\begin{matrix}
\xymatrix@!C@C-30pt{
\Ospin_5(M)    \ar[dr] \ar@/^2pc/[rr]   &&
\htildeM  \ar[dr] \ar@/^2pc/[rr] &&
\hept{B}\ar[dr]^{\beta}  \\
& \Ospin_5(B) \ar[dr] \ar[ur]  &&
\hept{M}  \ar[dr]^{\alpha}\ar[ur] &&\Ospin_4(B) \\
 \pi_1(\heqpt{B}) \ar[ur] \ar@/_2pc/[rr] && \rOspin_5(B,M)
  \ar[ur]^\gamma \ar@/_2pc/[rr]&&
\rOspin_4(M)\ar[ur]
}\end{matrix}
$$
\vskip .8cm\noindent
of exact sequences, valid
for any closed, oriented smooth or topological spin $4$-manifold $M$
(see Theorem \ref{mainnonspin} for 
the  analogous statement in the  non-spin case).
The maps labelled $\alpha$ and  $\beta$  are
not homomorphisms, so exactness  is understood in the
sense of  ``pointed sets" (meaning that $image = kernel$, where $kernel$ is the pre-image of the base point).

For the special case
when the fundamental group 
$\pi_1(M,x_0)$ is finite of odd order, we can compute 
this braid to obtain an explicit formula for $\hept{M}$
and a description of $\hM$.
The simply-connected case was already known
(see  \cite{ch1}, \cite{wallbook},\cite{lawson1},
\cite{kreck1}), but our proof is new even in that case.

In \cite[p.~85]{hk2} we defined the \emph{quadratic $2$-type}
of $M$ as the $4$-tuple $\quadtypeM$, where $s_M$ is the
intersection form on $\pi_2(M)$. The isometries 
of the quadratic $2$-type consist of all pairs of isomorphisms
 $\chi\colon\pi_1(M,x_0)\to \pi_1(M,x_0)$ and 
$\phi\colon\pi_2(M)\to \pi_2(M)$, such that $\phi(gx) = \chi(g)\phi(x)$,
which preserve the $k$-invariant and the intersection form.
The group $\Isom(\quadtypeM)$ is thus a  subgroup
of the arithmetic group $SO(\pi_2(M), s_M)$.

\begin{thma} Let $M^4$ be a connected, closed, oriented 
 smooth (or topological) manifold of dimension $4$. 
If $\pi_1(M,x_0)$ has odd
order, then
$$\hept{M} \cong \Kw \rtimes \Isom(\quadtypeM)$$
where $\Kw:=\ker(w_2\colon H_2(M;\cy 2) \to \cy
2))$.
\end{thma}
\noindent
The image of $\htildeM$ in $\hept{M}$ is isomorphic to
$\Isom(\quadtypeM)$, giving the semi-direct product splitting,
and the
 action  on the normal subgroup
$\Kw$ by $\Isom(\quadtypeM)$ is induced by the action of $\hept{M}$
on homology.

Let  $\cS^h(M\times I, \bd)$ denote the
 structure group of smooth or topological
manifold structures on $M\times I$, relative to
the given structure on $\bd(M\times I)$.
Let $\tilde L_6(\bZ[\pi])$ denote the \emph{reduced} Wall
group (see \cite[Chap.~9]{wallbook}) defined 
as the cokernel of the split injection
$L_6(\bZ) \to L_n(\bZ[\pi])$ induced by the inclusion
$1\to \pi$ of the trivial group. The group $\hM$ of smooth or topological  $h$-cobordisms from $M$ to $M$
is now determined up to extensions.
\begin{thmb} Let $M^4$ be a connected, closed,
oriented smooth (or topological)
 manifold of dimension $4$. If $\pi_1(M,x_0)$ has odd
order, then
there is a short exact sequence of groups:
$$1 \to \cS^h(M\times I,\bd) \to \hM \to 
\Isom(\quadtypeM) \to 1$$
where the normal subgroup $\cS^h(M\times I,\bd)$ is
abelian and is determined up to extension
by the
% (split) 
short exact sequence
$$0\to \tilde L_6(\bZ[\pi_1(M,x_0)]) \to \cS^h(M\times I,\bd)
\to H_1(M;\bZ) \to 0$$
of groups and homomorphisms.
\end{thmb} 
\noindent
In the simply-connected case, $ \cS^h(M\times I,\bd) =0$ so
the group of $h$-cobordisms is just isomorphic to the isometries
of the intersection form of $M$.

The authors would like to thank the referee for useful comments 
on the first version of this paper.

\section{Spin bordism groups and exact sequences}
We now define more precisely the objects and maps which appear in our braid,
beginning with the case when $M$ is a spin $4$-manifold.
We fix a lift  $\nu_M \colon M \to BSpin$ of the classifying map
for the  stable normal bundle of $M$.
Let $\Ospin_*(M)$ or $\Ospin_*(B)$ denote the singular bordism groups
 of topological
spin manifolds equipped with a reference map to $M$ or $B$.
The discussion below holds for smooth bordism 
(when $M$ is a smooth $4$-manifold) without  any essential
changes.

\subsection{The map $\alpha$.}
We fix a base-point $x_0 \in M$ and the corresponding base-point
in $B$, thinking of $B$ as constructed from $M$ by adding cells of
dimension $\geq 4$. The evaluation map at $x_0$ gives a fibration
$$\heqpt{M} \to \heq M \to M$$
where $\heq M$ denotes the space of orientation-preserving
homotopy self-equivalences of $M$.
We have a long exact sequence
$$\dots\to \pi_1(\heq{M})  \to \pi_1(M,x_0) \to \pi_0(\heqpt{M}) \to \pi_0(\heq{M})
\to \pi_0(M)\ .$$
The image $G(M,x_0):= \Image(ev_*\colon \pi_1(\heq{M})  \to \pi_1(M,x_0))$
has been studied by Gottlieb \cite{gottlieb1}. It is always a central
subgroup of $\pi_1(M,x_0)$, and 
$G(M,x_0)$ is trivial if $\chi(M) \neq 0$ (e.g. when
$\pi_1(M,x_0)$ is finite).

Since $M$ is connected, we see that any homotopy equivalence
is homotopic to a base-point preserving homotopy equivalence.
We are studying $\hept{M} := \pi_0(\heqpt{M})$. Notice
that the composition 
$$\pi_1(M,x_0) \to \pi_0(\heqpt{M}) \to \text{Aut}(\pi_1(M,x_0))$$
just sends an element $\sigma \in \pi_1(M,x_0)$ to the
automorphism ``conjugation by $\sigma$".

The inclusion gives a fixed reference map $c\colon
M \to B$ which is base-point preserving, and induces a homomorphism
$\hept{M} \to \hept{B}$, by obstruction theory.
We also have a  map
$$\alpha\colon \hept{M} \to \Ospin_4(M)$$
defined by $\alpha(f):= [M,f] -[M,id]$, but this is \emph{not} a
homomorphism:
\eqncount
\begin{equation}\label{alphamap}
 \alpha(f\circ g) = \alpha(f) + f_*(\alpha(g))  \ .
\end{equation}
Since $f$ is orientation-preserving, the fundamental
class $f_*[M]=[M]$ and so the image of
$\alpha$ is contained in the modified bordism group $\rOspin_4(M)$.

We have already mentioned the modified relative bordism groups 
$\rOspin_5(B,M)$, where the representing objects $(W,F)$ are spin manifolds
of dimension 5 with boundary, such that $f=F|\bd W$ has degree zero.
The usual bordism exact sequence of the pair $(B,M)$ can be adapted
to include these modified groups.
\begin{lemma}
There is an exact sequence 
$$\dots \Ospin_5(M) \to \Ospin_5(B) \to \rOspin_5(B,M) \to \rOspin_4(M)
\to \Ospin_4(B)\ .$$
\end{lemma}
\begin{proof} Left to the reader.
\end{proof}

\subsection{The groups $\htildeM$.}
Next we define the  groups $\htildeM$ as the  bordism groups
of objects $(W,F)$ where $W$ is a compact $5$-dimensional spin manifold
with $\bd_1 W = -M$ and $\bd_2 W = M$, and $F\colon W \to M$
is a continuous map such that $F|\bd_1 W = id_M$ and $F|\bd_2 W=f$
is a base-point and orientation-preserving homotopy equivalence. 
In particular, we mean that the spin structure on $W$ is a lift of $\nu_W$ to $BSpin$ which
agrees with our fixed lift for $\nu_M$ on both boundary components $\bd_1 W$ and $\bd_2 W$. We do not, however, require  the self-equivalence $f$ to preserve the spin structure on $M$.
Two such objects 
$(W,F)$ and $(W',F')$ are
\emph{bordant} if there is a base-point preserving homotopy $h$ between
$f=F|\bd_2  W$ and $f'=F'|\bd_2 W'$, such that  
the closed, spin $5$-manifold 
\eqncount
\begin{equation}\label{bordant}
 (-W' \cup_{\bd_1 W'=\bd_1 W} W\cup_{\bd_2 W =M\times 
0\disjointunion
\bd_2 W'=M\times 1} M\times I, F'\cup F\cup h)
\end{equation}
 represents zero in $\Ospin_5(M)$. We define a group
structure on $\htildeM$ by the formula
\eqncount
\begin{equation}\label{mult}
(W,F) \bullet (W',F') := (W\cup_{\bd_2 W=\bd_1 W' }
W', F \cup f\circ F')\ .
\end{equation}
 This is easily seen to be well-defined, and the inverse of
$(W,F)$ is represented by $(-W, f^{-1}\circ F)$ where $f^{-1}$ is a base-point
preserving homotopy inverse
for $f=F|\bd_2 W$. By convention, $\bd_1(-W) = \bd_2(W)$, so to
obtain an object of the form required we must adjoin a collar $M\times I$
to $-W$ along $\bd_1(-W)$
mapped into $M$ by a homotopy between $f^{-1}\circ f$ and $id_M$.
The different choices of such a homotopy result in bordant representatives
for the inverse.
The identity element in this group structure is
represented by the bordism $(M\times I, p_1)$, where $p_1\colon M\times I
\to M$ is the projection on the first factor. There is a homomorphism
$\Ospin_5(M) \to \htildeM$ by taking the disjoint union of a closed, spin
$5$-manifold mapping into $M$ and the identity element $(M\times I,p_1)$, and a homomorphism $\htildeM \to \hept{M}$ mapping
$(W, F)$ to the homotopy class of $f:=F|\bd_2 W$.
\begin{lemma}
There is an exact sequence of pointed sets
$$\xymatrix{\Ospin_5(M) \ar[r]& \htildeM \ar[r]& \hept{M} \ar[r]^{\alpha}&
\rOspin_4(M)}$$ 
where only the last map $\alpha$ fails to be a group homomorphism.
\end{lemma}
\begin{proof} Left to the reader.
\end{proof}

\subsection{The groups $\htildeB$.} 
To obtain a similar exact sequence through $\hept
B$ we start  by defining the group $\htildeB$ as the bordism group of objects
  $(W,F)$ where
$W$ is a compact
$5$-dimensional spin manifold with $\bd_1 W = -M$
 and $\bd_2 W = M$, and $F\colon
W \to B$ is a continuous map such that $F|\bd_1 W = c$ and $F|\bd_2 W=f$
is a base-point preserving $3$-equivalence. 
Two such objects 
$(W,F)$ and $(W',F')$ are
\emph{bordant} if there is a base-point preserving homotopy $h$ into $B$
between
$f=F|\bd_2  W$ and $f'=F'|\bd_2 W'$, such that  
the closed, spin $5$-manifold (\ref{bordant})
 represents zero in $\Ospin_5(B)$. To define the group structure on
$\htildeB$ we first remark that a base-point preserving $3$-equivalence
$f\colon M \to B$ induces a base-point preserving self-equivalence
$\phi_f\colon B \to B$ such that $\phi_f\circ c=f$. 
Furthermore, the map $\phi_f$ is uniquely defined by this equation, 
up to a base-point
preserving homotopy.

The multiplication   is now defined as in (\ref{mult}) by
the  formula 
$$(W,F) \bullet (W',F') := (W\cup_{\bd_2 W=
\bd_1 W' } W', F \cup \phi_f\circ F')$$
and the identity element  is represented by $(M\times I, c\circ p_1)$.
The inverse of
$(W,F)$ is represented by $(-W, \phi_f^{-1}\circ F)$ where $\phi_f^{-1}$ is a
base-point preserving homotopy inverse
for the self-equivalence $\phi_f\colon B\to B$
induced by $f=F|\bd_2 W$. 

\begin{lemma}\label{equality}
$\htildeM \cong \htildeB \ .$
\end{lemma}
\begin{proof}
We have a well-defined homomorphism $\htildeM \to \htildeB$
by composing with our reference map $c\colon M \to B$. Suppose
that $(W,F)$ represents an element in $\htildeB$.
Since $F|\bd_1 W = c$, the identity map on $M$ is a lift of
$F|\bd_1 W = c$ over $M$. We want to extend this lift over $W$
 by  homotoping the map $F\colon W \to B$ into
$M$, relative to $\bd_1 W$. By low-dimensional surgery on the map $F$, we may
assume that $F$ is $2$-connected.
 If $X$ denotes the homotopy fibre
of $c$, we are looking for the obstructions to lifting the map
$F$ relative to $\bd_1 W$ to the total space of the fibration
$X \to M \to B$. But the fibre $X$ is $2$-connected, so the
lifting obstructions lie in $H^{i+1}(W,\bd_1 W;\pi_i(X))$ for
$i \geq 3$. However, $H^{i+1}(W,\bd_1 W) = H_{5-i-1}(W;\bd_2 W) =0$
if $i \geq 3$, for any coefficients, since $F$ is $2$-connected and $f$ is
$3$-connected. If $\hat F\colon W \to M$ is a lift of $F$, then
$\hat f=\hat F|\bd_2 M$ is a 3-equivalence, and has degree 1
since it is bordant over $M$ to the identity map. Therefore $\hat f$ is
an orientation and base-point preserving homotopy equivalence.
This proves that the natural map $\htildeM \to \htildeB$
is \emph{surjective}. 

Suppose now that $(W,F)$ and $(W',F')$ represent two elements in
$\htildeM$ which are bordant over $B$. We may assume that the 
reference map $T\to B$ for the bordism
 is $3$-connected (by surgery on the interior of $T$), 
and then it follows as above that
there are no obstructions to lifting this reference map to $M$, relative to
the union of the boundary components $(W,F)$, $(W',F')$ and $\bd_1 W \times I$.
A  lifting of the reference map $T \to B$ restricted to $\bd_2 W \times I$
gives a base-point preserving homotopy between $f$ and $f'$
as required. Therefore $(W,F)$ and $(W',F')$  are bordant over $M$,
and the map $\htildeM \to \htildeB$ is \emph{injective}.
\end{proof}
\subsection{The map $\beta$.}
We have a map of pointed sets 
$$\beta\colon\hept B \to \Ospin_4(B)$$
defined by $\beta(\phi\colon B \to B) := [M,\phi\circ c] - [M,c]$. We
also have a homomorphism
$$\pi_1(\heqpt{B}) \to \Ospin_5(B)$$
 sending the adjoint map $h\colon B\times S^1 \to B$, for a representative
of an element in $\pi_1\heqpt{B}$,
to the bordism element $[M \times S^1, h\circ (c\times id)]$.
We use the null-bordant spin structure on the $S^1$ factor.
To see that this map induces a group
homomorphism, consider the surface $F$ obtained from the
$2$-disk by removing two small open balls in the interior. If
$h, h'\colon B\times S^1 \to B$ are the adjoints of maps
representing elements
of $\pi_1(\heqpt{B})$, and $h'' = h\bullet h'$ is the adjoint of the 
product, then there is an obvious map from $B \times F \to B$
such that the restriction to the boundary is given by $h$, $h'$
on the boundaries of the interior balls and by $h''$ on the exterior
boundary component. Then $M \times F$ gives the required spin bordism. 

Finally, the homomorphism
$\Ospin_5(M) \to \htildeM$ is defined taking the disjoint union of a closed, spin
$5$-manifold mapping into $B$ and the identity element $(M\times I,c\circ p_1)$.
\begin{lemma}
There is an exact sequence of pointed sets
$$\xymatrix{
\pi_1(\heqpt{B})\ar[r]& \Ospin_5(B) \ar[r]&
 \htildeM \ar[r]& \hept B \ar[r]^{\beta}& \Ospin_4(B)}$$
where only the last map $\beta$ fails to be a group homomorphism.
\end{lemma}
\begin{proof}
We will prove exactness for the related sequence
where $\htildeM$ is replaced by $\htildeB$, and then apply
Lemma \ref{equality}.
It follows easily from the definitions 
that the composite of any two
maps in this new sequence is trivial, and that we have 
 exactness at the terms $\htildeB$ and $\hept B$.
It remains to check exactness at $\Ospin_5(B)$. Let $(N,g)$ represent an
element of $\Ospin_5(B)$ which maps to the identity in $\htildeB$. 
This means that
the bordism element $(N,g) \disjointunion (M\times I,c\circ p_1)$
is bordant to $(M\times I, c\circ p_1)$. In particular, there is a base-point
preserving homotopy $h \colon M\times I \to B$ with $h|(M\times 0) =c$
and $h|(M\times 1)=c$. This homotopy  $h$ induces a pointed homotopy
$\hat h\colon B \times I \to B$ from the identity to the identity,
representing an element of $\pi_1(\heqpt{B})$.
\end{proof}
\subsection{The map $\gamma$.}
The remaining exact sequence in our braid diagram involves the
construction of a map $\gamma\colon \rOspin_5(B,M) \to \hept M$.
Let $(W,F)$ denote an element of $\rOspin_5(B,M)$. 
This is a $5$-dimensional spin manifold with boundary $(W,\bd W)$,
equipped with a reference map $F\colon W \to B$ such
that $F|\bd W$ factors through the classifying map $c\colon M \to B$.
We may assume that $\bd W$ is connected.

By taking the
boundary connected sum with the zero bordant element $(M\times I, p_1)$
along $\bd W$ and $M\times 1$, we may assume that
 $W$ has two boundary components $\bd_1 W =-M$ and $\bd_2 W = N$
with the reference map $F|\bd_1 W = c$. 
We may assume
(by low-dimensional surgery on the map $F$) that  $F$ 
is a $2$-equivalence. Moreover, since
$(W,F)$ is a modified bordism element, $g:=F|\bd_2 W$ is a degree 1
map from $N \to M$. 

Now consider the obstructions to lifting the
map $F\colon W \to B$ to $M$, relative to $F|\bd_2 W$. These
obstructions lie in the groups $H^{i+1}(W,N;\pi_i(X))$, where
$X$  denotes the fibre of the map $c\colon M \to B$. Since
$X$ is $2$-connected, after applying  Poincar\'e duality
we see that the lifting obstructions lie in
the  groups $H_{5-i-1}(W,M;\pi_i(X))=0$, for $i \geq 3$.
Let $r\colon W \to M$ be a 
lift of $F$ relative to $N$, and
consider the map $f:=r|\bd_1 W\colon M \to M$. Since $c\circ r \simeq F$,
and $F|\bd_1 W = c$, we have $c\circ f \simeq c$ and hence $f$
is a $3$-equivalence. However $f_*[M] = g_*[N]=[M]$ since $g$ has
degree 1 and the maps $f$ and $g$ are bordant over $M$.
However, a degree 1 map $f\colon M \to M$ which is a $3$-equivalence
is a homotopy equivalence (by Poincar\'e duality and Whitehead's Theorem).
We may assume that $f$ is also base-point preserving

\begin{lemma}\label{gamma}
There is a well-defined map
$$\gamma\colon \rOspin_5(B,M) \to \hept M\ .$$
\end{lemma}
\begin{proof}
Let $r\colon W \to M$ be a lifting of $F$ and let
$f:=r|\bd_1 W$.
We define
$$\gamma(W,F) := [f\colon M \to M] \in \hept M\ .$$
To see that the map $\gamma$ is well-defined, suppose that $(W',F')$
is another representative for the same relative bordism class and that
we have already found liftings $r$ and $r'$ of the maps $F$ and $F'$
respectively. Let $(T, \varphi)$
denote a bordism between $(W,F)$ and $(W', F')$, respecting
the boundary. More
precisely, $\bd T$ consists of the union of $W$, $W'$, $M\times I$ and
a $5$-dimensional bordism $P$ between $N$ and $N'$.
We may assume that the reference map $\varphi\colon T \to B$
 is a $3$-equivalence
by surgery on the interior of $T$. Now consider the obstructions to lifting
the map $\varphi$ to $M$, relative both to $\varphi | P$ and to
our chosen liftings $r$ and $r'$. 
The obstructions to lifting $\varphi$ lie in the groups
$H^{i+1}(T, W\cup P\cup W';\pi_i(X))$ for $i \geq 3$. 
They may be evaluated by
Poincar\'e duality as above, and are again zero. Then any such lifting
$\hat\varphi\colon T \to M$ of $\varphi$ gives a homotopy 
$h=\hat\varphi |(M\times I)$ between
$f=r|\bd_1 W$ and $f'=r'|\bd_1 W'$. 
We may assume in addition that $h$ preserves
the base-point $x_0$,  by constructing our lifting relative to a thickening
$D^4\times I\subset M \times I$ of the interval $x_0 \times I$.
\end{proof}
To check that the map $\gamma$ fits into our braid diagram, we introduce
another object. Let $\hB$ denote the equivalence classes of triples
$(M\times I,h,f)$, where $f\colon M \to M$ is a base-point
preserving homotopy equivalence, and
 $h\colon M\times I
\to B$ is a base-point preserving homotopy 
between $c$ and $c\circ f$.  Two triples 
$(M\times I, h,f)$ and $(M\times I, h', f')$ are equivalent
if there is a base-point preserving 
homotopy $p\colon M \times I \to M$ between
$f$ and $f'$, and a continuous map $t\colon M\times I \times I \to B$
such that $t|M\times I \times 0 = h$, $t|M\times I \times 1 = h'$,
$t|M\times 0\times I = c$ and $t|M\times 1 \times I = c\circ p$.
We define a multiplication on $\hB$ by the union
$$ (M\times I, h,f)\bullet (M\times I, h',f') = 
(M\times I, h \cup h'\circ f, f'\negthinspace\circ\negthinspace f)$$
where the two copies of $M \times I$ on the left-hand side 
are identified with $M \times
[0,1/2]$ and $M \times [1/2,1]$ respectively on the right-hand side. 
The inverse of
$(M\times I, h, f)$ is represented by
$(M\times I, \bar{h}\circ f^{-1},f^{-1})$, where
$\bar{h}(x,t) = h(x, 1-t)$, for $0\leq t \leq 1$ and $x\in M$, and $f^{-1}$ is a
base-point preserving homotopy inverse for
$f$. We adjoin a pointed homotopy
 between $f\circ f^{-1}$ at the end $M \times 0$
to obtain an element in the standard form.
\begin{lemma}\label{Bexact}
There is an exact sequence of groups and homomorphisms
$$\xymatrix{\pi_1(\heqpt B) \ar[r]& \hB \ar[r]^{\bd}& 
\hept M \ar[r]& \hept B}\ .$$
\end{lemma}
\begin{proof}
It is clear that the map $(M\times I, h,f)\mapsto f$ gives a homomorphism
$\bd\colon\hB \to \hept M$, and the exactness is just a formal consequence of the
definitions. In particular, $\Image \bd$ is a normal subgroup of
$\hept{M}$.
\end{proof}
\begin{remark}\label{conjugation}
 For each $[g] \in \hept{M}$,  we have a 
base-point preserving self-equivalence
$\phi_g\colon B \to B$ such that $c\circ g= \phi_g\circ c$.
There is a conjugation action  on $\hB$ defined by
$$(M\times I, h,f)\mapsto 
(M\times I, \phi_g\circ h\circ g^{-1}, g\circ f\circ g^{-1} )$$
which is compatible with the boundary map $\bd\colon \hB \to
\hept{M}$. 
\end{remark}
\begin{lemma}\label{eta}
 There is a bijection
$\eta\colon\rOspin_5(B,M)\cong \hB $ such that $\bd\circ\eta=\gamma$.
\end{lemma}
\begin{proof}
Let $(W,F)$ represent an element of $\rOspin_5(B,M)$,
as constructed at the beginning of this sub-section,  with $F$ a 2-connected
map as usual. We have constructed a lifting $r\colon W \to M$ of $F$
relative to $g:=F|\bd_2 W$. 
In other words, $F$ is homotopic to $r$ over $B$
and we can use a homotopy to give a map
 $\varphi\colon W \times I \to B$
such that $\varphi| W\times 0 = F$,
$\varphi | W\times 1 = r$, and $\varphi|\bd_2\times I = g$. Let 
$h:=\varphi |\bd_1 W \times I$. Then $h\colon M\times I$ is a homotopy
between $c$ and $c\circ f$ where $f:=r|\bd_1 W$.
We define a map $\eta\colon\rOspin_5(B,M)\to \hB $ by
 $(W,F) \mapsto (M\times I, h, f)$. It is easy to check that this map is well-defined and 
gives a bijection between the two
sets. By construction, the map $\gamma=\bd\circ\eta$ is the composite of this bijection 
and the boundary map $\bd\colon\hB \to \hept M$ from
Lemma \ref{Bexact}. 
\end{proof}
\begin{remark}
This argument shows that the element $(W,F)$ is bordant 
to $(M\times I, h)$,
so it represents the same bordism class in $\rOspin_5(B,M)$. However, we
do not know if the bordism group structure (addition by disjoint union) agrees
with the multiplication $\bullet$ defined on the elements
$(M\times I, h)$. For this reason, we don't know if the map $\gamma$
is always a homomorphism. If $\pi_1(M,x_0)$ has odd order,
it turns out that $\gamma$ \emph{is} a homomorphism.
\end{remark}
\begin{corollary}
There is an exact sequence of pointed sets
$$\xymatrix{
\pi_1(\heqpt B) \ar[r] &\rOspin_5(B,M) \ar[r]^{\gamma}&
\hept M \ar[r] &\hept B}\ .$$
\end{corollary}
\begin{proof}
Left to the reader.
\end{proof}
\subsection{Commutativity of the braid} We have now verified the exactness
of all the sequences in the braid diagram, so it remains to check the
commutativity of the diagram. We will only discuss two of the sub-diagrams.
\begin{lemma}
The composite $\alpha\circ \gamma$ equals the boundary
map $\bd\colon\rOspin_5(B,M) \to \rOspin_4(M)$. 
\end{lemma}
\begin{proof}
Let $(W,F)$ represent an element of $\rOspin_5(B,M) $ in the standard
form above. Then its image in $\rOspin_4(M)$ is represented
by $[N,g] - [M,id]$, where $g:=F|\bd_2 W$ as usual. However, the
existence of a lifting $r\colon W \to M$ for $F$ shows that
$[N ,g]$ is bordant over $M$ to $[M,f]$, and so $\bd(W,F)$ represents the same bordism element as
$\alpha\circ\gamma(W,F)$.
\end{proof}
\begin{lemma}
The composite $\Ospin_5(B) \to \htildeM \to \hept M$
equals the composite  $\Ospin_5(B) \to \rOspin_5(B,M) \to \hept M$ 
up to inversion $[f]\mapsto [f]^{-1}$ in $\hept M$.
\end{lemma}
\begin{proof}
Let $(N,g)$ denote an element of $\Ospin_5(B)$. Then we map it
into $\htildeM$ by forming the connected sum 
$(M\times I \Sharp N, p_1\Sharp g)$
and lifting the map $ \varphi=p_1\Sharp g$ to 
$\hat\varphi\colon M\times I \Sharp N \to M$ relative to
$M\times 0$. Then $[\hat\varphi | M\times 1]$ is the image of the first
composition in $\hept M$.
To compute the other composition, we again form the connected sum
$(M\times I \Sharp N, p_1\Sharp g)$ and lift the map 
$ \varphi=p_1\Sharp g$ to $r\colon M\times I \Sharp N \to M$ relative
to $M\times 1$. Then the image of the second composition
is represented  by $f:=r|M\times 0$. However, notice that the map
$f^{-1}\circ r$, together with a pointed homotopy from $f^{-1}\circ f$
in a small collar of $M \times 0$, gives another lifting of $\varphi$ relative
to $M \times 0$. Therefore $[\hat\varphi | M\times 1] = [f^{-1}]$,
showing that the two compositions agree up to inversion in $\hept M$.
\end{proof}
We have proved that our braid diagram is
\emph{sign-commutative}, meaning that the sub-diagrams are all strictly
commutative except for the two composites ending in $\hept M$ which
only agree up to inversion. 
\begin{theorem}\label{mainspin}
Let $M$ be a closed, oriented smooth (respectively topological)
$4$-manifold. If $M$ is a spin manifold, there is a sign-commutative
diagram of exact
sequences
\vskip .2cm
$$\begin{matrix}
\xymatrix@!C@C-30pt{
\Ospin_5(M)    \ar[dr] \ar@/^2pc/[rr]   &&
\htildeM  \ar[dr] \ar@/^2pc/[rr] &&
\hept{B}\ar[dr]^{\beta}  \\
& \Ospin_5(B) \ar[dr] \ar[ur]  &&
\hept{M}  \ar[dr]^{\alpha}\ar[ur] &&\Ospin_4(B) \\
 \pi_1(\heqpt{B}) \ar[ur] \ar@/_2pc/[rr] && \rOspin_5(B,M)
  \ar[ur]^\gamma \ar@/_2pc/[rr]&&
\rOspin_4(M)\ar[ur]
}\end{matrix}
$$
\vskip .8cm\noindent
involving the bordism groups of smooth (respectively topological)
spin manifolds. All the maps except $\alpha$, $\beta$ and 
possibly $\gamma$
are group homomorphisms.
\end{theorem}
\section{Non-spin bordism groups}
When $w_2(M) \neq 0$ the bordism groups must be modified in
the above braid in order to carry out the  arguments used
to establish commutativity.

Let $\xi\colon E \to BSO$ be a fibration, and recall that elements in the bordism groups
$\Omega_n(E)$ are represented by maps 
$\bar \nu\colon N \to E$ from  a smooth, closed,  $n$-manifold $N$
into $E$, such that
 $\xi \circ \bar \nu = \nu_N$, where $\nu_N \colon N \to BSO$
classifies the stable normal bundle of $N$. The bordism relation
also involves a compatible lifting of the normal bundle data over
the cobordism (see \cite{lashof1}, \cite[Chap.~II]{stong1}).

Recall that a \emph{normal $k$-smoothing} of $M$ in $E$ is
a lifting $\bar\nu\colon M \to E$ of $\nu_M$  such that $\bar\nu$
is a $(k+1)$-equivalence \cite[p.~711]{kreck3}. The fibration
$E\to BSO$ is called $k$-universal if  its fibre is connected,
with homotopy groups vanishing in dimensions $\geq k+1$.
The \emph{normal $2$-type} of $M$ is a $2$-universal
fibration $E \to BSO$ admitting a normal $2$-smoothing
of $M$ (see \cite[p.~711]{kreck3} for an extensive development of these concepts).

For the non-spin case of our braid we will use the bordism
groups of the normal $2$-type:

$$\xymatrix@!C
{BSpin\ar[r]^i\ar@{=}[d]&
\Bw\ar[r]^j\ar[d]^{\xi} &
B\ar[d]^{w_2}\\ 
BSpin\ar[r] & BSO\ar[r]^{w} & \ K(\cy 2,2)
}$$
as described in \cite[\S 2]{teichner1}.  The map $w = w_2(\gamma)$
pulls back the second Stiefel-Whitney class for the universal
oriented vector bundle
$\gamma$ over $BSO$. The
``James" spectral sequence used to compute $\Omega_*(\Bw)
=\pi_*(M\xi)$
has the same $E_2$-term as the one used above for
$w_2=0$, but the differentials are twisted by $w_2$.
In particular, $d_2$ is the dual of $Sq^2_w$,
where $Sq^2_w(x):= Sq^2(x) + x\cup w_2$.
 There is a corresponding non-spin version of $\Ospin_*(M)$,
namely the bordism groups $\Omega_*(\Mw):=\pi_*(M\xi)$
of the Thom space
associated to the fibration:
$$\xymatrix@!C
{BSpin\ar[r]^i\ar@{=}[d]&
\Mw\ar[r]^j\ar[d]^{\xi} &
M\ar[d]^{w_2}\\ 
BSpin\ar[r] & BSO\ar[r]^{w} & \ K(\cy 2,2)
}$$
Again the $E_2$-term of the James spectral sequence
is unchanged from the spin case, but the differentials are twisted  by
$w_2$ with the above formula for $Sq^2_w$.
As in the spin case, we choose a particular representative for
the map $w_2$ such that $w_2 = w\circ \nu_M$.

Our next step is to define a suitable ``thickening" of $\hept{M}$
for the non-spin case. Here is the main technical ingredient.
\begin{lemma}\label{liftone}
Let $f\colon M \to M$ be a base-point and
orientation-preserving homotopy equivalence.
Then there exists a base-point preserving homotopy equivalence
$f'\colon M \to M$, such that
 $f\simeq f'$ preserving the base point, with
 $w\circ \nu_M = w_2\circ f'$.
\end{lemma}
\begin{proof}
By the Dold-Whitney Theorem \cite{doldwhitney1}, there
is an isomorphism $f^*(\nu_M) \cong \nu_M$. We therefore
have a (base-point preserving) homotopy $h\colon M \times I \to BSO$
between the classifying maps  $\nu_M \circ f \simeq \nu_M$. 
Now define
$\hat f\colon M \to \Mw$ lifting $\nu_M\circ f$ by the formula
$$\hat f (x) := (f(x), \nu_M(f(x))$$
for all $x\in M$, and note that this makes sense because $w_2= w\circ
\nu_M$ as maps to $K(\cy 2,2)$.
We apply the covering homotopy theorem to get
$\hat h \colon M \times I \to \Mw$
lifting $h$, with the property that
 $\xi \circ (\hat h\vv_{M\times 1}) = \nu_M$.
Let $f'\colon M \to M$ be defined by the formula
$f':=j\circ (h\vv_{M\times 1})$, where $j\colon \Mw \to M$
is the projection on the first factor.
Then $f'\simeq f$ by the homotopy $j\circ \hat h$, and
 we have
$$w_2(f'(x)) = w_2(j(\hat h\vv_{M\times 1}(x))) = 
w(\xi(\hat h\vv_{M\times 1}(x))) = w(\nu_M(x))$$
for all $x\in M$, as required.
\end{proof}
As a consequence of  the Lemma,
  the formula $\hat f' (x):= (f'(x), \nu_M(x))$
gives a map $\hat f'\colon M \to \Mw$, and in fact $\hat f' \equiv
 \hat h\vv_{M\times 1}$ by construction.
Therefore $\xi\circ \hat f' = \nu_M$ as maps $M \to BSO$.
We will consider the set of all such maps into $\Mw$ under
a suitable equivalence relation.

\begin{definition}\label{thickaut}
Let $\whept{M}$ denote the set of 
equivalence classes of maps $\hat f\colon M \to \Mw$ such that
(i) $f:=j\circ \hat f $ is a base-point and orientation preserving
homotopy equivalence, and (ii) $\xi\circ \hat f =\nu_M$. Two such maps
$\hat f$ and $\hat g$ are \emph{equivalent} if there exists a homotopy
$\hat h\colon M\times I \to \Mw$ such that $h:=j\circ \hat h$ is
a base-point preserving homotopy between $f$ and $g$, and
$\xi\circ \hat h =\nu_M\circ p_1$, where $p_1\colon M\times I \to M$
denotes projection on the first factor.
\end{definition}
Given two maps $\hat f, \hat g\colon M \to \Mw$ as above,  we define
$$\hat f \bullet \hat g\colon M \to \Mw$$
as the unique map from $M$ into the pull-back $\Mw$ defined  by the pair
$f\circ g\colon M \to M$ and $\nu_M\colon M \to BSO$. Since
$w_2\circ f\circ g = w\circ \nu_M\circ g = w_2\circ g = w\circ \nu_M$,
this pair of maps is compatible with the pull-back.
\begin{lemma} $\whept{M}$ is a group under this operation.
\end{lemma}
\begin{proof} To check that the operation just defined passes to
equivalence classes, suppose that $\hat h$ is a homotopy as above
between $\hat f$ and $\hat f'$ representing the same element of
$\whept{M}$. Let $h:=j\circ \hat h$ and notice that $w_2\circ h\circ g
= w\circ \nu_M \circ p_1\circ (g\times id) = w_2\circ g\circ p_1 = w\circ \nu_M\circ p_1$. We have a similar argument in the case when
$\hat g$ is varied by a homotopy.

Next we discuss the identity element and inverses. 
Let $\widehat{id}_M\colon M \to \Mw$ denote the map defined
by the pair $(id_M\colon M \to M, \nu_M\colon M \to BSO)$. This
map will represent the identity element in our group structure.

Given $\hat f$ representing an element
of $\whept{M}$, let $\hat g\colon M \to \Mw$ be a map
constructed as in Lemma \ref{liftone} applied to any base-point
preserving homotopy inverse $f^{-1}$ for $f:=j\circ \hat f$.
Now if $h\colon M \times I \to M$ is a base-point preserving
homotopy between $f\circ g$ and $id_M$, we can assume
that $w_2\circ h = w_2\circ p_1$.
To see this, note that the different
maps $M\times I \to K(\cy 2,2)$ relative to the given maps on the boundary are classified by $H^1(M;\cy 2)$. But we can construct
a map $M\times S^1 \to M$ using any element of $\pi_1(M,x_0)$,
and this gives a homotopy from $id_M$ to itself realizing any
desired element of $H^1(M;\cy 2)$.  It follows that the pair of maps
$h\colon M \times I \to M$ and $\nu_M\circ p_1\colon M\times I
\to BSO$ define a unique map $\hat h \colon M\times I \to \Mw$.
This is exactly the required homotopy between $\hat f\circ\hat g$
and $\widehat{id}_M$. We will refer to $\hat h$ as an
\emph{admissible} homotopy. Checking the remaining properties of the group
structure will be left to the reader.
\end{proof}

Now we will define a map
$$\alpha\colon \whept{M} \to  \widehat\Omega_4(\Mw)$$
for use in our braid, where the modified bordism groups are defined
by letting the \emph{degree} of a reference map $\hat g\colon N^4 \to \Mw$ be the ordinary degree of $g:=j\circ \hat g$.
Given $[\hat f] \in \whept{M}$, let
$$\alpha(\hat f) := [M, \hat f] - [M, \widehat{id}_M] \in \Omega_4(\Mw),$$
and notice that this element has degree zero.
Since $\xi\circ \hat f = \nu_M$, we have a bundle map
$\hat b\colon \nu_M \to \xi$ and a commutative diagram
$$\xymatrix{ E(\nu_M) \ar[r]^{\hat b}\ar[d] & E( \xi)\ar[d]\cr
M \ar[r]^{\hat f} & \Mw}
$$
expressing that fact that $(M,\hat f)$ represents an element of
the bordism theory for the normal 2-type. It is clear from
the way that the equivalence relation is defined for $\whept{M}$
that $\alpha$ is well-defined, independent of the choice 
of representative for $[\hat f]$.

Next comes the definition of $\whtildeM$ and the homomorphism
$\whtildeM \to \whept{M}$.

\begin{definition} Let $\whtildeM$ denote the bordism groups of
 pairs $(W, \widehat F)$, where $W$ is a compact,
oriented $5$-manifold with $\bd_1 W = -M$ and $\bd_2 W = M$.
The map $\widehat F \colon W \to \Mw$ restricts to $\widehat{id}_M$
on $\bd_1 W$, and on $\bd_2 W$ to a map $\hat f\colon M \to \Mw$ satisfying
properties (i) and (ii) of Definition \ref{thickaut} .
\end{definition}
Two such objects 
$(W,\widehat F)$ and $(W',\widehat F')$ are
\emph{bordant} if there is an equivalence $\hat h$ between
$\hat f=\widehat F|\bd_2  W$ and $\hat f'=\widehat F'|\bd_2 W'$, such that  
the closed $5$-manifold 
\eqncount
\begin{equation}\label{wbordant}
 (-W' \cup_{\bd_1 W'=\bd_1 W} W\cup_{\bd_2 W =M\times 
0\disjointunion
\bd_2 W'=M\times 1} M\times I, \widehat F'\cup \widehat F\cup h)
\end{equation}
 represents zero in $\Omega_5(\Mw)$. We define a group
structure on $\whtildeM$ by the formula
\eqncount
\begin{equation}\label{wmult}
(W,\widehat F) \bullet (W',\widehat F') := (W\cup_{\bd_2 W=\bd_1 W' }
W', \widehat F \cup \hat f\bullet \widehat F')\ .
\end{equation}
 This is easily seen to be well-defined, and the inverse of
$(W, \widehat F)$ is represented by $(-W, \hat f^{-1}\bullet\widehat F)$ where $\hat f^{-1}$ represents the
inverse
for $\hat f=\widehat F|\bd_2 W$ in $\whept{M}$.
 By convention, $\bd_1(-W) = \bd_2(W)$, so to
obtain an object of the form required we must adjoin a collar $M\times I$
to $-W$ along $\bd_1(-W)$
mapped into $M$ by an admissible homotopy between $\hat f^{-1}\bullet \hat f$ and $\widehat{id}_M$.
The different choices of such a homotopy result in bordant representatives
for the inverse.
The identity element in this group structure is
represented by the bordism $(M\times I, \hat p_1)$, where $\hat p_1:=\widehat{id}_M \circ p_1$ and
$p_1\colon M\times I
\to M$ is the projection on the first factor. There is a homomorphism
$\Omega_5(\Mw) \to \whtildeM$ by taking the disjoint union of a closed, 
$5$-manifold with normal structure in $\Mw$ and the identity element $(M\times I,\hat p_1)$.
\begin{lemma}
There is an exact sequence of pointed sets
$$\xymatrix{\Omega_5(\Mw) \ar[r]& \whtildeM \ar[r]& \whept M \ar[r]^{\alpha}&
\widehat\Omega_4(\Mw)}$$ 
where only the last map $\alpha$ fails to be a group homomorphism.
\end{lemma}
\begin{proof}
The homomorphism $\whtildeM \to \whept{M}$
is defined on representatives by sending $(W, \widehat F)$ to
$\hat f: = \widehat F|\bd_2 W$. The rest of the details will be left to the reader.
\end{proof} 

Finally, we will define the analogous bordism groups $\whtildeB$
and the group $\whept{B}$ of self-equivalences, together with the map
$\beta\colon \whept{B} \to  \Omega_4(\Bw)$.
Here is the basic technical ingredient.
\begin{lemma}\label{techB}
Given a base-point preserving map $f\colon M \to B$, there is a unique
extension (up to base-point preserving homotopy)
 $\phi_f\colon B \to B$ such that
$\phi_f\circ c = f$. If $f$ is a $3$-equivalence then $\phi_f$
is a homotopy equivalence. If $w_2\circ f = w_2$, then
$w_2\circ \phi_f = w_2$.
\end{lemma}
\begin{proof}
The existence and uniqueness of the extension $\phi_f$ follow
from obstruction theory, since $c\colon M \to B$ is a $3$-equivalence.
The other statements are clear.
\end{proof}
\begin{definition}\label{thickautB}
Let $\whept{B}$ denote the set of 
equivalence classes of maps $\hat f\colon M \to \Bw$ such that
(i) $f:=j\circ \hat f $ is a base-point preserving
$3$-equivalence, and (ii) $\xi\circ \hat f =\nu_M$. Two such maps
$\hat f$ and $\hat g$ are \emph{equivalent} if there exists a homotopy
$\hat h\colon M\times I \to \Bw$ such that $h:=j\circ \hat h$ is
a base-point preserving homotopy between $f$ and $g$, and
$\xi\circ \hat h =\nu_M\circ p_1$, where $p_1\colon M\times I \to M$
denotes projection on the first factor.
\end{definition}
Given two maps $\hat f, \hat g\colon M \to \Bw$ as above,  we define
$$\hat f \bullet \hat g\colon M \to \Bw$$
as the unique map from $M$ into the pull-back $\Bw$ defined  by the pair
$\phi_f\circ \phi_g\circ c\colon M \to B$
 and $\nu_M\colon M \to BSO$. Here we are using Lemma \ref{techB} 
to factor the maps $\phi_f\circ c = f$ and $\phi_g\circ c = g$.
Since $w_2\circ \phi_f\circ \phi_g \circ c
= w\circ \nu_M\circ \phi_g \circ c
= w_2\circ \phi_g\circ c = w\circ \nu_M$,
this pair of maps is compatible with the pull-back.
\begin{lemma} $\whept{B}$ is a group under this operation.
\end{lemma}
\begin{proof}
Let $\hat{c}\colon M \to \Bw$ denote the map defined
by the pair $(c\colon M \to B, \nu_M\colon M \to BSO)$. This
map will represent the identity element in our group structure.

Given $\hat f$ representing an element
of $\whept{B}$, write $f=\phi_f\circ c$ as above and choose
a base-point preserving homotopy inverse $\psi\colon B \to B$
for $\phi_f$, with the additional property that $w_2\circ \psi = w_2$.
This is another ``lifting" argument using the fibration $B\to K(\cy 2,2)$.
Then the pair $g:= \psi\circ c$ and $\nu_M$ define a map
$\hat g\colon M \to \Bw$ representing the inverse of $\hat f$.
We leave the check that $\hat g \bullet \hat f\simeq \hat c $ via an admissible equivalence to the reader.
\end{proof}

\begin{definition} Let $\whtildeB$ denote the bordism groups of
 pairs $(W, \widehat F)$, where $W$ is a compact,
oriented $5$-manifold with $\bd_1 W = -M$ and $\bd_2 W = M$.
The map $\widehat F \colon W \to \Bw$ restricts to $\hat{c}$
on $\bd_1 W$, and on $\bd_2 W$ to a map $\hat f\colon M \to \Bw$ satisfying
properties (i) and (ii) of Definition \ref{thickautB} .
\end{definition}
Two such objects 
$(W,\widehat F)$ and $(W',\widehat F')$ are
\emph{bordant} if there is an equivalence $\hat h$ between
$\hat f=\widehat F|\bd_2  W$ and $\hat f'=\widehat F'|\bd_2 W'$, such that  
the closed $5$-manifold (\ref{wbordant})
 represents zero in $\Omega_5(\Bw)$. We define a group
structure on $\whtildeB$ as in (\ref{wmult}) by the formula
\eqncount
\begin{equation}
(W,\widehat F) \bullet (W',\widehat F') := (W\cup_{\bd_2 W=\bd_1 W' }
W', \widehat F \cup \hat f\bullet \widehat F')\ .
\end{equation}
and the identity element is represented by $(M\times I, \hat  p_1)$,
where $\hat p_1:= \hat c\circ p_1$.
 The inverse of
$(W, \widehat F)$ is represented by $(-W, \hat f^{-1}\bullet\widehat F)$ where $\hat f^{-1}$ represents the
inverse
for $\hat f=\widehat F|\bd_2 W$ in $\whept{B}$.
\begin{lemma}\label{equalitynonspin}
 $\whtildeM \cong \whtildeB$.
\end{lemma}
\begin{proof}
This  follows as  in the proof of Lemma \ref{equality} for the spin case:
the lifting arguments take place over the fixed map $\nu_M\colon
M\to BSO$.
\end{proof}

There is a homomorphism
$\Omega_5(\Bw) \to \whtildeB$ by taking the disjoint union of a closed, 
$5$-manifold with normal structure in $\Bw$ and the identity element $(M\times I,\hat p_1)$. Furthermore, we have a map
$\beta\colon \whept{B} \to \Omega_4(\Bw)$ defined by
$\beta(\hat f) := [M, \hat f] - [M, \hat c]$.

We can also define $\heqpt{M,w_2}$ and $\heqpt{B,w_2}$ as the spaces of maps from $M \to \Mw$ or $M \to \Bw$
satisfying the properties (i) and (ii) of Definitions \ref{thickaut} or \ref{thickautB} respectively. Then $\whept{M} = \pi_0(\heqpt{M,w_2})$ and $\whept{B} = \pi_0(\heqpt{B,w_2})$. We therefore have a homomorphism
$\pi_1(\heqpt{B,w_2}) \to \Omega_5(\Bw)$ sending the adjoint map
$\hat h\colon M \times S^1 \to \Bw$ for a representative of an element
in $\pi_1(\heqpt{B,w_2})$ to the bordism element 
$(M\times S^1, \hat h)$ in the normal $2$-type $\Bw$.

\begin{lemma} There is an exact sequence of pointed sets
$$\xymatrix@C-5pt{\pi_1(\heqpt{B,w_2})\ar[r] & \Omega_5(\Bw) \ar[r]& \whtildeM 
\ar[r]&
 \whept{B} \ar[r]^{\beta}& \Omega_4(\Bw)}$$
where only the last map $\beta$ fails to be a group homomorphism.
\end{lemma}
\begin{proof} Left to the reader.
\end{proof}

These definitions and properties allow is to establish our commutative
braid. 
\begin{theorem}\label{mainnonspin}
Let $M$ be a closed, oriented smooth (respectively topological)
$4$-manifold with normal $2$-type $\Bw$. There is a sign-commutative
diagram of exact
sequences
\vskip .2cm
$$\begin{matrix}
\xymatrix@C-30pt{
\Omega_5(\Mw)    \ar[dr] \ar@/^2pc/[rr]   &&
\whtildeM  \ar[dr] \ar@/^2pc/[rr] &&
\whept{B}\ar[dr]^{\beta}  \\
& \Omega_5(\Bw) \ar[dr] \ar[ur]  &&
\whept{M}  \ar[dr]^{\alpha}\ar[ur] &&\Omega_4(\Bw) \\
 \pi_1(\heqpt{B,w_2}) \ar[ur] \ar@/_2pc/[rr] && \widehat\Omega_5(\Bw,\Mw)
  \ar[ur]^\gamma \ar@/_2pc/[rr]&&
\widehat\Omega_4(\Mw)\ar[ur]
}\end{matrix}
$$
\vskip .8cm
\noindent
involving the bordism groups of smooth (respectively topological)
manifolds. 
\end{theorem}
As before,
the two composites ending in $\whept{M}$ agree up to inversion,
and the other sub-diagrams are strictly commutative.
\begin{proof} The proof of this result follows the pattern for the
spin case. The key points are Lemma \ref{equalitynonspin} and
the definition of the map $\gamma$ (see Lemma \ref{gamma}).
\end{proof}
We conclude this section by pointing out the connection between
$\hept{M}$ and $\whept{M}$.
\begin{lemma} There is a short exact sequence of groups
$$0\to H^1(M;\cy 2) \to \whept{M} \to \hept{M} \to 1\ .$$
\end{lemma}
\begin{proof}
There is a natural map $\heqpt{M,w_2} \to \heqpt{M}$ defined by sending $\hat f$ to $f:=j\circ \hat f$, and this induces a
surjective homomorphism on the groups of homotopy classes over $BSO$.
The identification of the kernel with $H^1(M;\cy 2)$ follows from
the fibration $K(\cy 2,1) \to \Mw \to M\times BSO$ and
obstruction theory.
\end{proof}
\begin{remark}
A similar result holds for $\whept{B}$, which maps surjectively
onto the subgroup of $\hept{B}$ fixing $w_2$. The kernel
is again isomorphic to $H^1(M;\cy 2)$.
\end{remark}

\section{Odd order fundamental groups}
In this section, we assume 
 that $\pi_1(M,x_0)$ is a finite group of odd order.
We can then compute the terms in our braid to obtain a
more explicit expression for $\hept{M}\cong \whept{M}$. 
We also have some information about the group $\hM$ of $h$-cobordisms.

Notice that (even without assumption on $\pi_1(M,x_0)$) there is an exact sequence
$$1\to \cS^h(M\times I, \bd) \to \hM \to \hept{M}$$
where $\cS^h(M\times I, \bd)$ denotes the
 structure group of smooth or topological
manifold structures on $M\times I$, relative to
the given structure on $\bd(M\times I)$.

We first point out a useful input from surgery theory.
\begin{lemma}\label{Hone}
Suppose that $\pi_1(M)$ is finite of odd order.
 There is an injection
$H_1(M;\bZ) \to \whtildeM$, factoring through the map
$\Omega_5(\Mw) \to \whtildeM$ from the braid diagram.
\end{lemma}
\begin{proof}
 We have a commutative diagram of exact sequences
$$\xymatrix{&\tilde L_6(\bZ[\pi_1])\ar@{=}[r] \ar[d]&\tilde L_6(\bZ[\pi_1])\ar[d]&\cr
0\ar[r]&\cS^h(M\times I,\bd)\ar[r]\ar[d]&\hM\ar[r]\ar[d]&\hept{M}\cr
&H_1(M;\bZ)\ar[r]\ar[d]&\whtildeM\ar[d]&\cr
&L_5(\bZ[\pi_1]) \ar@{=}[r]&L_5(\bZ[\pi_1]) &\cr
}
$$
\smallskip\noindent
where the left-hand vertical sequence is from Wall's surgery
exact sequence \cite[Chap.~10]{wallbook}. To obtain the
right-hand vertical sequence we use the modified surgery
theory of \cite{kreck3}. The surgery
obstruction map $\whtildeM \to L_5(\bZ[\pi_1])$ from \cite[Thm.~4]{kreck3} is the obstruction to finding a bordism
over the normal type to an element of $\hM$ (see the remark on \cite[p.~734]{kreck3}
to replace the monoid $\ell_5(\bZ[\pi_1])$ by the Wall group,
and the remark on  \cite[p.~738]{kreck3} for the $h$-cobordism
version). By construction, this map is a homomorphism. There is also an action of $L_6(\bZ[\pi_1])$ on 
$\hM$, as in the surgery exact sequence, which again by
construction gives a homomorphism. The exactness of 
the displayed right-hand sequence  follows from
\cite[Thm.~3]{kreck3} and the remark \cite[p.~730]{kreck3}.

The horizontal maps come from the bordism interpretation of
the surgery exact sequence
$$L_6(\bZ[\pi_1]) \to \cS^h(M\times I,\bd)
\to \cT(M\times I, \bd) \to L_5(\bZ[\pi_1])$$
in which the normal invariant term $\cT(M\times I, \bd)$ is the set
of degree 1 normal maps $F\colon (W,\bd W) \to (M\times I, \bd)$, inducing
the identity on the boundary \cite[Prop.~10.2]{wallbook}.
The group structure on this set is defined as for 
$\whtildeM$.
The map  $\cT(M\times I, \bd) \to \whtildeM$ takes such an element to $(W, \widehat F)\in \whtildeM$. This map factors through  $\Omega_5(\Mw)$ by sending such an element to the bordism class of $(W\cup M\times I, \widehat F)$.
On the other hand, there is an isomorphism of groups
$$\cT(M\times I, \bd) \cong [M\times I, \bd;G/TOP]
 = [SM;G/TOP]$$
when we use the $co$-$H$-space
structure on the reduced suspension 
$SM$ of $M$.
 That group structure agrees
with the usual one for the normal invariants from the $H$-space structure on $G/TOP$ (see \cite[\S 1.6]{spanier}).

 A computation gives $[M\times I, \bd;G/TOP] \cong H_1(M;\bZ)$, 
and a diagram
chase now shows that the composite map 
$H_1(M;\bZ)\to \whtildeM$ is an injection.
\end{proof}

\begin{remark} For later use, we will note that the map
$H_1(M;\bZ) \to \Omega_5(\Mw)$ defined above may be identified with the  homomorphism 
$$H_1(M;\bZ) = E_2^{1,4} \to E_\infty^{1,4} \subset
\Omega_5(\Mw)$$
 in the Atiyah-Hirzebruch spectral
sequence whose $E_2$-term is $H_p(M;\Ospin_q(\ast))$.

To see this, we consider an embedding $f\colon (S^1 \times D^3)\times I \to M\times I$ representing an element of $H_1(M;\bZ)$. There is a commutative diagram
$$\xymatrix{\cT(S^1\times D^4,\bd)\ar[r]\ar[d] &
\Ospin_5(S^1)\ar[d]\cr
\cT(M\times I,\bd)\ar[r]&
\Omega_5(\Mw)
}$$
where the left vertical map is given by gluing a normal map with
range $S^1 \times D^4 \equiv S^1\times D^3\times I$ into 
$M\times I$, and extending by the identity.
By Poincar\'e duality, there is a
commutative diagram 
$$\xymatrix{\cT(S^1\times D^4,\bd)\ar[d]&
H^4(S^1\times D^4, \bd)\ar[l]_{\approx}\ar[r]^{\approx}\ar[d]^{f^!}&
H_1(S^1;\bZ)\ar[r]\ar[d]^{f_*}&\Ospin_5(S^1)\ar[d]\cr
\cT(M\times I,\bd)&H^4(M\times I, \bd)\ar[l]_{\approx}\ar[r]^{\approx}&
H_1(M;\bZ)\ar[r]&\Omega_5(\Mw)
}
$$
factoring the one above,
where $f^!$ denotes the map induced by the collapse
$M\times I \to S^1\times D^4/S^1\times S^3$.
The identification of $H_1(M;\bZ)$ with the normal invariants uses  Poincar\'e duality with $L$-spectrum coefficients, but in this low
dimensional situation it reduces to the ordinary duality.
The last horizontal maps in this diagram are induced from
the maps $E_2^{1,4} \to E_\infty^{1,4}$ in the spectral sequences.
\end{remark}

The remaining proofs will be done in a number of steps, 
starting with the case
of \emph{spin} manifolds.
We mean \emph{topological} bordism throughout and
homology with integral coefficients unless otherwise noted.
\begin{proposition}
Let $B$ denote the normal $2$-type of a spin $4$-manifold
$M$ with odd order fundamental group. Then
$\Ospin_4(B) \subset H_4(B) \oplus \bZ$ and
there is a short exact sequence
$0 \to  H_1(M) 
\to \Ospin_5(B)\to H_5(B)$.
\end{proposition}
\begin{proof}
This follows from the Atiyah-Hirzebruch spectral sequence,
whose $E_2$-term is $H_p(B;\Ospin_q(\ast))$.
The first differential 
$d_2\colon E_2^{p,q} \to E_2^{p-2,q+1}$ is given
by the dual of $Sq^2$ (if $q=1$) or this composed with reduction
mod 2 (if $q=0$), see \cite[p.~751]{teichner1}. 
We substitute the values
$\Ospin_q(\ast) = \bZ, \cy{2}, \cy{2}, 0, \bZ, 0$,
for $0 \leq q \leq 5$. Then the differential for $(p,q)= (4,1)$ becomes
$d_2\colon H_4(B;\cy 2) \to H_2(B; \cy 2)$.
This homomorphism may be detected by transfer to
the universal covering $\tilde B$, since $\pi_1$ has odd
order. Notice that $\tilde B$ is just a product of $\bC P^\infty$'s.
It follows that 
$Sq^2\colon H^2(\tilde B;\cy 2) \to H^4(\tilde B;\cy 2)$
is injective, hence its dual is surjective even when restricted
to the subgroup of $\pi_1$-invariant elements (by averaging).
Therefore, on the line $p+q=4$, the only groups which
survive to $E_\infty$ are $\bZ$ in the $(0,4)$ position,
and a subgroup of $H_4(B)$ in the $(4,0)$ position. 

For the line $p+q=5$, we have again that the differential
$d_2\colon H_6(B;\bZ) \to H_4(B; \cy 2)$ from position
$(p,q) = (6,0)$ is surjective
onto the kernel of the above differential
 $d_2\colon H_4(B;\cy 2) \to H_2(B; \cy 2)$. 
This follows from the exactness of the sequence
$$\xymatrix{H^2(\tilde B;\cy 2)\ar[r]^{Sq^2}& 
H^4(\tilde B;\cy 2)\ar[r]^{Sq^2}&H^6(\tilde B;\cy 2)}$$ 
and the surjectivity of $H_6(B;\bZ) \to H_6(B;\cy 2)$.
Finally, by transfer to $\tilde B$ we get
$H_3(B;\cy 2) =0$. Therefore the groups that survive on this
line are
$H_1(B) = H_1(M)$ in the $(1,4)$ position (by Lemma \ref{Hone})
and $H_5(B)$ in 
the $(5,0)$ position. 
\end{proof}
\begin{lemma}
$\ker (\beta\colon \hept{B} \to \Ospin_4(B)) \subseteq \Isom
([\pi_1,\pi_2, k,s])$.
\end{lemma}
\begin{proof}
Although $\beta$ is not a homomorphism, we can still
define $\ker(\beta) = \beta^{-1}(0)$.
The natural map $\Ospin_4(B) \to H_4(B)$ sends
a bordism element to the image of its fundamental class.
If $\phi \in \hept{B}$, and $c\colon M \to B$ is its classifying
map, then $\beta(\phi):= [M, \phi\circ c] - [M, c]$.
The image of this element in $H_4(B)$ is zero when
$\phi_*(c_*[M]) =c_*[M]$. But $\trf(c_*[M]) = s(M)$, 
the intersection form of $M$ on $\pi_2$ considered
as an element in $H_4(B) = \Gamma(\pi_2)$,  so 
$\ker\beta$ is contained in the self-equivalences of $B$ which
preserve the quadratic $2$-type. 
\end{proof}
Next we calculate some more bordism groups and determine the
image of the map
$\alpha\colon \hept{M} \to \rOspin_4(M)$.
\begin{proposition}
$\Ospin_4(M) = \bZ \oplus H_2(M;\cy 2)\oplus \bZ$, and
$\Image\alpha = H_2(M;\cy 2)$.
$\Ospin_5(M)=\cy 2\oplus H_1(M)$. 
The map $\Ospin_5(M)
\to
\Ospin_5(B)$ is projection onto the subgroup $H_1(M)$.
\end{proposition}
\begin{proof}
We use the same spectral sequence, but the terms are
a bit simpler because $H_3(M)=H_3(M;\cy 2)=H_5(M)=0$,
and $H_4(M;\cy 2) =\cy 2$. 
Since $M$ is spin, the map
$Sq^2\colon H^2(M;\cy 2) \to H^4(M;\cy 2)$ is zero,
so the differential $d_2\colon E_2^{4,1} \to E_2^{2,2}$
is also zero.
 The  line
$p+q=4$ now gives 
$\Ospin_4(M) = \bZ \oplus H_2(M;\cy 2) \oplus H_4(M) $.
If $f\colon M \to M$ represents an element of
$\hept{M}$, then $\alpha(f):= [M,f] - [M, id]$. It follows
that $\alpha(f) \in H_2(M;\cy 2)$ since both the signature
and the fundamental class in $H_4(M)$ are preserved
by a homotopy equivalence.

For the line $p+q=5$ in the $E_2$-term, we have $H_1(M)$ in
the
$(1,4)$ position, and $H_4(M;\cy 2)\cong \cy 2$ in the 
$(4,1)$ position. 
and both these terms survive to $E_\infty$. Under the
map $\Ospin_5(M) \to \Ospin_5(B)$, the summand
$H_1(M)$ maps isomorphically (by Lemma \ref{Hone} again), 
and the $\cy 2$
summand
maps to zero.
It follows that $H_2(M;\cy 2)$ lies in the image from
$\rOspin_5(B,M) \to \rOspin_4(M)$, hence $\Image \alpha =
H_2(M;\cy 2)$.
\end{proof}
\begin{corollary}
$\ker(\htildeM \to \hept{B}) = \ker(\htildeM \to \hept{M}) 
\cong H_1(M)$.
\end{corollary}

Now we need to compute some homology groups.
We need the following special case of a result of P. Teichner.
\begin{lemma}[\cite{teichner3}]
If $M$ has odd order fundamental group, then
$\Gamma(\pi_2(M)) = \bZ\oplus
\pi_3(M)$ as $\Lambda:=\bZ[\pi_1(M)]$ modules.
\end{lemma}
\begin{proof}
Recall that $s_M\in \Gamma(\pi_2)$ denotes the equivariant intersection
form on $\pi_2(M)$. The $\pi_1$-module $\Gamma(\pi_2)$ sits in the
Whitehead sequence
$$0 \to H_4(\widetilde M;\bZ) \to \Gamma(\pi_2) \to \pi_3(M) \to 0$$
where the first map sends $1\in H_4(\widetilde M;\bZ)$
to $s_M$ (see \cite{hk2}).
We wish to construct a $\pi_1$-homomorphism
$f\colon  \Gamma(\pi_2) \to \bZ$ such that $f(s_M) = 1$.
First, consider the map $f_1\colon \Gamma(\pi_2)\to \bZ$ given by 
the composite of the norm $N\colon \Gamma(\pi_2)
\to H^0(\pi_1;\Gamma(\pi_2))$ and a map $d\colon
H^0(\pi_1;\Gamma(\pi_2)) \to \bZ$ chosen so that
$d(s_M) =1$. Such a map exists because $s_M$ is unimodular
and is thus a primitive element in $\Gamma(\pi_2)$. We have
$f_1(s_M) = |\pi_1|$.

Next, a map $f_2\colon \Gamma(\pi_2) \to \bZ$ was constructed
in \cite{bauer1}, as the composite
$$\Gamma(\pi_2)\subseteq \pi_2\otimes\pi_2 \cong
\Hom(\pi_2^*, \pi_2) \cong \Hom(\pi_2,\pi_2) \to \bZ$$
where the middle isomorphisms are defined by
$x\otimes y \mapsto (\psi\mapsto \psi(x)\cdot y)$ and
$\theta(\psi):= \psi\circ s_M^{-1}$,  and the last map
is the trace. By definition, $f_2(s_M) = trace(id_{\pi_2}) =
\rank \pi_2$. But $\rank \pi_2 = \chi(\widetilde M) - 2 = |\pi_1|\cdot\chi(M) -2$
is relatively prime to $|\pi_1|$, so we can get a homomorphism
$f\colon \Gamma(\pi_2) \to \bZ$ with $f(s_M)=1$ by taking an
appropriate linear combination of $f_1$ and $f_2$.
\end{proof}
\begin{proposition}
$H_4(B)$ is torsion-free, and $H_5(B) =0$.
\end{proposition}
\begin{proof}
We use the Serre spectral sequence of the fibration
$\tilde B \to B \to K(\pi,1)$, 
 $E_2$-term
given by $E_2^{p,q} = H_p(\pi_1;H_q(\tilde B))$,
where
$\pi_1=\pi_1(M)$, and substitute the values
$H_i(\tilde B)=0$, for $i=1,3,5$, $H_2(\tilde B)
=\pi_2:=\pi_2(M)$, and $H_4(\tilde B)=\Gamma(\pi_2)$.
We have a splitting $\Gamma(\pi_2) = \bZ\oplus
\pi_3(M)$ as $\Lambda:=\bZ[\pi_1]$ modules. But from \cite[p.~3]{bauer1} we have
$Tors(H_4(B)) \cong \widehat H_0(\pi_1, \pi_3(M))$.
Also, from  \cite[\S 3]{hk6}, and the
assumption that $\pi_1$ has odd order, we have
$\widehat H^i(\pi_1;\Gamma(\pi_2)) = \widehat
H^i(\pi_1;\bZ)$ in all dimensions. In particular,
$\widehat H_0(\pi_1, \pi_3(M))=0$ implying that
$H_4(B)$ is torsion-free, and
the term $E_2^{1,4} = H_1(\pi_1;\Gamma(\pi_2))
= H_1(\pi_1;\bZ)$.

The image
of the projection map $H_4(\tilde B) \to H_4(B)$ is always a quotient of $E_2^{0,4}=H_0(\pi_1;H_4(\tilde B))$ under
the edge homomorphism. 
In our case,  
$\widehat H^{-1}(\pi_1;\bZ)=0$ (for any finite group
\cite{brown1}), so we have an inclusion
$$H_0(\pi_1;\Gamma(\pi_2)) \subset
H^0(\pi_1;\Gamma(\pi_2))\ .$$
 But $\Gamma(\pi_2)$
is $\bZ$-torsion free, hence so is the term
$E_2^{0,4}=H_0(\pi_1;H_4(\tilde B))$.
It follows that this term survives to $E_\infty$ and
injects into $H_4(B)$.

Now consider the differentials $d_3$ in the spectral sequence
affecting the lines $p+q=4,5$. These have the form
$d_3\colon H_{i+3}(\pi_1) \to H_{i}(\pi_1;\pi_2)$,\
for $i=1,2$ or $3$.
We can obtain information about them  by comparing the spectral
sequence for $B$ with that for $B_2$, the $2$-skeleton
of $B$ in some $CW$-structure. By the results of 
\cite[\S 2]{hk2}, we have $\pi_2(B_2) = \Omega^3\bZ$,
and a short exact sequence of stable $\Lambda$-modules
$$0\to \Omega^3\bZ \to \pi_2 \to S^3\bZ \to 0\ .$$
However, the corresponding differentials in the spectral
sequence for $B_2$ must be isomorphisms  (in order
that $H_*(B_2) =0$ for $*>2$). We can therefore identity our
original $d_3$ differentials with the natural maps
$$H_i(\pi_1;\Omega^3\bZ)
\to H_i(\pi_1;\pi_2)$$
 in the long exact sequence
$$ \dots\to H_{i+1}(\pi_1; S^3\bZ) \to
H_i(\pi_1;\Omega^3\bZ)
\to H_i(\pi_1;\pi_2) \to H_{i}(\pi_1; S^3\bZ)\to  \dots$$
 for the extension describing
$\pi_2$,
by means of the dimension-shifting isomorphism
$H_i(\pi_1;\Omega^3\bZ) = H_{i+3}(\pi_1)$.
Now we compute the maps in this long exact sequence,
using the values $H_3(\pi_1;S^3\bZ) = \hat H_0(\pi_1)
=0$, and $H_2(\pi_1;S^3\bZ)=\hat H_{-1}(\pi_1)
=\cy{|\pi_1|}$. Since $H_3(B)=0$ (following from
the fact that $H_3(M)=0$ and the $3$-equivalence
$M \to B$), a comparison with the $3$-skeleton $B_3$
shows that the differential
$d_3\colon  H_4(\pi_1) \to H_1(\pi_1;\pi_2)$
is an isomorphism.
We also 
get the following 
exact sequences:
$$H_4(\pi_1;\pi_2) \to  H_1(\pi_1) \to H_6(\pi_1) \to
H_3(\pi_1;\pi_2)\to 0$$ and 
$$0 \to H_5(\pi_1) \to H_2(\pi_1;\pi_2)\to\cy{|\pi_1|}
\to 0$$
determining the other $d_3$ differentials. 

Finally, by comparing to the spectral sequence for the
$4$-skeleton $B_4\subset B$ and the spectral sequence for the
universal covering $\widetilde M \to M$, we can see that the
differential
$d_3\colon E_2^{4,2} \to E_2^{1,4}$ is just the
natural map $H_4(\pi_1,\pi_2) \to H_1(\pi_1)$ above.
Furthermore, we can identify
the differential $d_5\colon E_3^{6,0} \to E_3^{1,4}$
with the inclusion $\ker(H_6(\pi_1) \to
H_3(\pi_1;\pi_2)) \subseteq H_1(\pi_1)$ given by
the exact sequence above. This eliminates everything
on the line $p+q=5$, so $H_5(B)=0$.
\end{proof}

\begin{corollary}\label{gamma_injects}
The group $\rOspin_5(B,M) = H_2(M;\cy 2)$ and
injects into $\hept{M}$.
\end{corollary}
\begin{proof}
We had a short exact sequence
$0\to H_1(M) \to \Ospin_5(B) \to H_5(B)$,
but now we know that $H_5(B)=0$. Therefore
$$\rOspin_5(B,M) = \ker(\rOspin_4(M)\to \Ospin_4(B)),$$
which equals $H_2(M;\cy 2)$. The result now
follows by the commutativity of the braid.
\end{proof}

\begin{corollary}
The images  of $\hept{M}$ or 
$\htildeM$ in $\hept{B}$
are precisely equal to the
isometry group  
$\Isom(\quadtypeM))$ of the quadratic $2$-type.
\end{corollary}
\begin{proof}
If $f\colon M \to M$ is an element in $\hept{M}$, then
its image in $\Ospin_4(B)$ factors through the map
$\rOspin_4(M) \to \Ospin_4(B)$, which has
trivial image in $H_4(B)$. Therefore, $c_*(f_*[M])
=c_*[M]$, and since $\trf(c_*[M])$ is just the
intersection  form of $M$ (considered as an
element of $H_4(\tilde B)$ \cite[p.~89]{hk2}), we see
that $\Image (\hept{M} \to \hept{B})$ is contained in
the isometries of the quadratic $2$-type.

However, since $H_4(B)$ is torsion free, it is detected by
the transfer map $\trf\colon H_4(B) \to H_4(\tilde B)$.
Now suppose that $\phi\colon B\to B$ is an element of
$\hept{B}$ contained in $\Isom(\quadtypeM))$.
Then 
$$\trf(\phi_*(c_*[M]))=\trf(c_*[M]),$$ and hence
$\phi_*(c_*[M])=c_*[M]$. By \cite[1.3]{hk2}, there exists
a lifting $h\colon M\to M$ such that $c\circ h\simeq \phi\circ
\phi$.  It follows (as in \cite[p.~88]{hk2}) that
$h$ is a homotopy equivalence.
The result for the image of $\htildeM$ follows
by  exactness of the braid, and the fact 
that $H_4(B)$ is torsion free.
\end{proof}
We can now put the pieces together to establish
our main results. Here are the relevant terms of our braid diagram:

\vskip .2cm
$$
\xymatrix@!C@C-54pt{
H_1(M)\oplus \cy{2}    \ar[dr] \ar@/^2pc/[rr]   &&
\htildeM  \ar[dr] \ar@/^2pc/[rr] &&
\Isom([\pi_1,\pi_2,k,s])\ar[dr]  \\
& H_1(M) \ar[dr]_{0}\ar[ur]  &&
\hept{M}  \ar[dr]^{\alpha}\ar[ur] &&0 \\
0 \ar[ur] \ar@/_2pc/[rr] && H_2(M;\cy 2)
  \ar[ur]^{\gamma} \ar@/_2pc/[rr]^{\cong}&&
H_2(M;\cy 2)\ar[ur]
}
$$
\bigskip
\begin{proof}[The proof of Theorem B]
We work in the topological category, and explain the smooth case in Remark \ref{toptosmooth}.
The first exact sequence (for spin manifolds)
$$1\to \cS(M\times I,\bd)\to \hM \to 
\Isom(\quadtypeM)\to 1$$
is obtained from the diagram in the proof of Lemma \ref{Hone}, by replacing
$\hept{M}$ with the image of $\htildeM$ in $\hept{M}$. From the braid diagram,
we see that this image is just the isometry group 
$\Isom(\quadtypeM)$.
The exact sequence for $\cS(M\times I,\bd)$ is
 part of the surgery exact sequence
\cite{wallbook}. We have just
substituted the calculation $L^h_5(\bZ[\pi_1(M)]) =0$
(see \cite{bak1}),
and
computed the normal invariant term 
$$[M\times I,M\times \bd I;G/TOP]
=H^2(M\times I,\bd ; \cy 2) \oplus H^4(M\times I,\bd;\bZ)\ .$$
The fact that $\tilde L_6(\bZ[\pi_1(M,x_0)])$
injects into  $\cS^h(M\times I,\bd)$ for odd order fundamental
groups is a computation of the surgery obstruction map
$$[M\times D^2,M\times S^1;G/TOP] \to L_6(\bZ[\pi_1(M,x_0)])$$
in the surgery exact sequence. This map factors through
a bordism group depending functorially on $\pi_1(M,x_0)$
(see \cite[Thm.~13B.3]{wallbook}). Since the $2$-localization map
$L_*(\bZ[\pi]) \to L_*(\bZ[\pi])\otimes \bZ_{(2)}$ is an injection for $L$-groups of finite groups \cite[Thm.~7.4]{wall-V}, we can use the fact that  $2$-local bordism 
is generated by the image from the $2$-Sylow subgroup.
 It follows that the image of the surgery obstruction
map $[M \times D^2,\bd;G/TOP] \to L_6(\bZ[\pi_1(M,x_0)])$
 factors through the $2$-Sylow subgroup inclusion $L_6(\bZ) \to L_6(\bZ[\pi_1(M,x_0)])$. We have given a direct argument here, but
this fact about the surgery obstruction map also follows from
\cite[Thm.~A]{hmtw1}.

In the non-spin case, we must still prove that the image
of $\htildeM$ in $\hept{B}$ is still $\Isom(\quadtypeM)$.
This will be done below.
\end{proof}
\begin{proof}[The proof of Theorem A]
In the spin case, 
the quotient of $\htildeM$ by the subgroup
$H_1(M)$ is isomorphic to $\Isom(\quadtypeM)$.
This gives the  splitting of the short exact
sequence 
$$0\to K \to \hept{M}
 \to \Isom(\quadtypeM)\to 1$$
where  $K:=\ker(\hept{M} \to \hept{B})$. It follows
that 
$$\hept{M} \cong K \rtimes \Isom(\quadtypeM)$$
with the conjugation action of $\Isom(\quadtypeM)$ on
the normal subgroup $K$ defining the semi-direct product
structure.
However, the braid diagram also shows that
the map $\gamma$ is an injective \emph{homomorphism}.
To check this, first observe that the isomorphism $\rOspin_4(M)= H_2(M;\cy 2)\oplus \bZ$ is natural, so any self-homotopy equivalence of $M$ which acts as the identity
 on $H_2(M;\cy 2)$ also acts as the identity on $\rOspin_4(M)$.
But any element in the image of $\gamma$ is trivial in $\hept{B}$, so acts as the identity on $H_2(M;\cy 2)$. Then 
formula (\ref{alphamap}) shows that $\alpha$ is a homomorphism
on the image of $\gamma$, and a diagram chase using Corollary \ref{gamma_injects} shows that $\gamma$ is a homomorphism.

 Therefore we  have a short
exact sequence of groups and homomorphisms
$$0\to H_2(M;\cy 2) \to \hept{M} \to \Isom(\quadtypeM)\to 1\ .$$
Moreover, 
$K = \Image\gamma$  and $K$ is mapped isomorphically onto
$H_2(M;\cy 2)$ by the map $\alpha$.
Finally, we apply formula (\ref{alphamap}) to obtain the
relations:
$$0 = \alpha(id_M) = \alpha(g\circ g^{-1}) =
\alpha(g) + g_*(\alpha(g^{-1}))$$
for any $[g]\in \hept{M}$, and
$$\alpha(g\circ f\circ g^{-1}) = g_*(\alpha(f))$$
for any $[f]\in K$. Therefore the conjugation action on $K$
agrees with the induced action on homology under 
the identification $K \cong H_2(M;\cy 2)$ via $\alpha$.
It follows that 
$$\hept{M}\cong H_2(M;\cy 2) \rtimes 
  \Isom(\quadtypeM)$$
as required, with the action of $\Isom(\quadtypeM)$ on the
normal subgroup $H_2(M;\cy 2)$ given by 
the induced action
of homotopy self-equivalences on homology. 
This completes the proof in the spin case.

\smallskip
For the non-spin case we must compute the bordism
groups of the normal $2$-type. Recall
that the first differential in the ``James" spectral
sequence used to compute $\Omega_*(\Bw)
=\pi_*(M\xi)$
has the same $E_2$-term as the one used above for
$w_2=0$, but the differentials are twisted by $w_2$.
In particular, $d_2$ is the dual of $Sq^2_w$,
where $Sq^2_w(x):= Sq^2(x) + x\cup w_2$.

\begin{proposition}
$\Omega_4(\Bw) = \bZ\oplus \cy 2\oplus
H_4(B)$ and $\Omega_5(\Bw)=H_1(M)$.
$\Omega_4(\Mw) = \bZ\oplus H_2(M;\cy 2)
\oplus \bZ$, and $\Omega_5(\Mw) 
=H_1(M)\oplus \cy 2$.
The natural map $\Omega_4(\Mw)
\to \Omega_4(\Bw)$
is injective on the $\bZ$ summands, and is
the homomorphism
$w_2\colon H_2(M;\cy 2) \to \cy 2$ on $H_2(M;\cy 2)$.
\end{proposition}
\begin{proof}
As before, we only need to compute the $d_2$
differentials. The point is that the composition
$$\xymatrix{H^2(B;\cy 2) \ar[r]^{Sq^2_w}&
H^4(B;\cy 2) \ar[r]^{Sq^2_w}&
H^6(B;\cy 2) }$$
is exact and the kernel of $Sq^2_w\colon 
H^2(B;\cy 2)\to H^4(B;\cy 2)$ is the subspace
$\langle w_2\rangle\cong \cy 2$. This gives the
cokernel $\cy 2$ in the $E^{2,2}_\infty$ position.
The same calculation in the spectral sequence for
$\Mw$ uses the fact
that 
$$Sq^2_w\colon 
H^2(M;\cy 2)\to H^4(M;\cy 2)$$
 is zero, since $w_2$
is also the first Wu class of $M$.
\end{proof}
We now continue with the proof of Theorem A
and Theorem B
in the non-spin case. 
 The
relevant terms on our braid are now:

\vskip .2cm
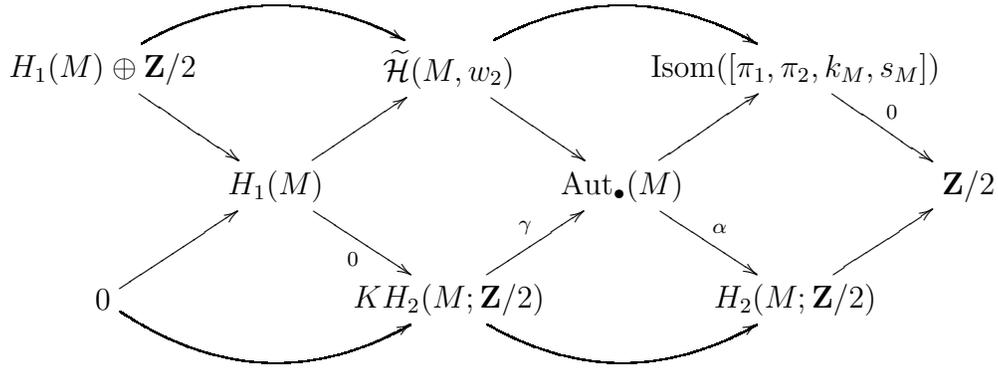
$$
\xymatrix@!C@C-75pt{
H_1(M)\oplus \cy{2}    \ar[dr] \ar@/^2pc/[rr]   &&
\whtildeM  \ar[dr] \ar@/^2pc/[rr] &&
\Isom(\quadtypeM)
\ar[dr]^{0}  \\
& H_1(M) \ar[dr]_{0}\ar[ur]  &&
\hept{M}  \ar[dr]^{\alpha}\ar[ur] &&\cy 2\\
0 \ar[ur] \ar@/_2pc/[rr] && \Kw
  \ar[ur]^{\gamma} \ar@/_2pc/[rr]&&
H_2(M;\cy 2)\ar[ur]
}
$$
\vskip .8cm
\noindent
Since the class $w_2 \in H^2(M;\cy 2)$
is a characteristic element for the cup product form
(mod 2), it is preserved by the induced map of a self-homotopy
equivalence of $M$. Therefore, the image of $\hept{M}$
in $\Omega_4(\Mw)$ lies in the
subgroup $\Kw :=\ker(w_2\colon H^2(M;\cy 2)\to
\cy 2)$. It then follows from the braid diagram, that
$$\Image(\hept{M} \to \hept{B}) =\Isom(\quadtypeM)$$
just as in the spin case. This
completes the proof of Theorem B, and Theorem A follows
as in the spin case. 
\end{proof}
\begin{remark}\label{toptosmooth}
If $M$ and $M'$ are
smooth, closed, oriented $4$-manifolds, and
$W$ is a topological $h$-cobordism between them, 
then there is a single obstruction in $H^4(W, \bd W;\cy 2)$
 to smoothing $W$ relative to the boundary
(see \cite[p.~194, 202]{kirby-siebenmann1}, or \cite[8.3B]{fq1}). If $\pi_1(M, x_0)$ has odd order
this obstruction vanishes. This implies that the forgetful map
$\cH_{\text{DIFF}}(M) \to \cH_{\text{TOP}}(M)$ is surjective.
It is also injective: we compare the smooth and topological
 surgery exact sequences
for $\cS^h(M\times I,\bd)$ as in Lemma \ref{Hone},
noting that the map on the term $H_1(M;\bZ)$ is multiplication by
$2$ (hence an isomorphism).
It follows that the  calculation in Theorem B also holds for the smooth $h$-cobordism group $\cH_{\text{DIFF}}(M)$. 
In addition, if $M$ and $M'$ are smooth $4$-manifolds
which are homeomorphic, 
then there exists a smooth $h$-cobordism
between them.
 It follows that the set $\cH (M, M')$ of \emph{smooth}
$h$-cobordisms between $M$ and $M'$ is in bijection with
$\cH_{\text{DIFF}}(M)$, whenever $\cH (M, M')$ is non-empty. In particular, 
$\cH (M, M')$ is also computed by Theorem B (extending the result of \cite{kreck1} for the simply-connected case).
\end{remark}
\bibliographystyle{amsplain}
\providecommand{\bysame}{\leavevmode\hbox to3em{\hrulefill}\thinspace}
\providecommand{\MR}{\relax\ifhmode\unskip\space\fi MR }
\providecommand{\MRhref}[2]{%
  \href{http://www.ams.org/mathscinet-getitem?mr=#1}{#2}
}
\providecommand{\href}[2]{#2}

\end{document}